\documentclass[a4paper,10pt]{article}
\usepackage{amssymb,amsmath}

\newcommand{\bal}{\mathbf{\alpha}}
\newcommand{\al}{\alpha}
\newcommand{\N}{\mathbb{N}}
\newcommand{\R}{\mathbb{R}}
\newcommand{\ep}{\varepsilon}
\newcommand{\beql}[1]{\begin{equation}\label{#1}}
\newcommand{\eeq}{\end{equation}}
\newcommand{\modd}[1]{\; ( \text{mod} \; #1)}
\newtheorem{theorem}{Theorem}
\newtheorem{lemma}{Lemma}

\begin{document}
\title{The Cubic Case of Vinogradov's Mean Value Theorem --- 
A Simplified Approach to Wooley's ``Efficient Congruencing''}
  
\author{D.R. Heath-Brown\\Mathematical Institute, Oxford}
\date{}
\maketitle
\section{Introduction}

In a remarkable series of papers, Wooley \cite{WA}, \cite{WE2},
\cite{WME}, \cite{WAp}, \cite{WC}, and in collaboration with
Ford  \cite{FW}, has made dramatic progress with Vinogradov's
mean value theorem. This has culminated very recently in the full
proof of the main conjecture, by Bourgain, Demeter and Guth
\cite{BDG}, using rather different methods.  
Wooley's survey article \cite{WS} gives an excellent
introduction to his results and their applications. The mean value theorem 
concerns the integer $J_{s,k}(X)$ defined as the number of solutions
$(x_1,\ldots,x_{2s})\in \N^{2s}$ of the simultaneous equations
\beql{eqs}
x_1^j+\ldots+x_s^j=x_{s+1}^j+\ldots+x_{2s}^j\;\;\;(1\le j\le k)
\eeq
with $x_1,\ldots,x_{2s}\le X$. Here $X\ge 1$ is an arbitrary real
number, and $s$ and $k$ are positive integers,
which one treats as being fixed.  The key feature of this system is
that if $(x_1,\ldots,x_{2s})$ is a solution, so is any translate
$(x_1+c,\ldots,x_{2s}+c)$.

The various forms of the Vinogradov mean value theorem give upper
bounds for $J_{s,k}(X)$. It is not hard to see that
\[J_{s,k}(X)\gg_{s,k} X^s+X^{2s-k(k+1)/2},\]
for $X\ge 1$, and the central conjecture is that
\[J_{s,k}(X)\ll_{s,k,\ep} X^{\ep}(X^s+X^{2s-k(k+1)/2})\]
for any $\ep>0$. ``Classically'' this was known for $k=1$ and $2$, for
$s\le k+1$, and for $s\ge s_0(k)$ with a value $s_0(k)\ll k^2\log k$.  However
Wooley \cite{WAp} shows that one may take $s_0(k)=k^2-k+1$, and that
the conjecture also holds for $s\le s_1(k)$ with
$s_1(k)=k(k+1)/2-k/3+o(k)$. Finally, in \cite{WC}, he shows that the
conjecture holds for $k=3$.

The purpose of this paper is to present a much simplified version of
Wooley's methods, sufficient to handle the case $k=3$. 
\begin{theorem}
We have
\[J_{6,3}(X)\ll_{\ep}X^{3+\ep}\]
for any fixed $\ep>0$.
\end{theorem}
It is trivial from (\ref{1}) below that if $s$ and $t$ are any
positive integers then we will have $J_{s+t,k}(X)\le
X^{2t}J_{s,k}(X)$ and $J_{s,k}(X)\le J_{s+t,k}(X)^{s/(s+t)}$.
Thus for
$k=3$ we can deduce the general case of the conjecture immediately
from the theorem.

It should be stressed that, while the argument of the present paper
appears cleaner than that presented by Wooley \cite{WC}, it is merely
a simplification of his version.  The underlying principles are the
same.  It is not a ``different'' proof. In part the simplification
arises from the restriction to the case $k=3$.  

\section{Outline of the Proof}

Investigations into the mean value theorem depend crucially on an
alternative interpretation of $J_{s,k}(X)$ in terms of exponential
sums. If $\bal\in\R^k$ we write
\[f_k(\bal;X)=f(\bal)=\sum_{x\le X}e(\al_1 x+\ldots+\al_k x^k),\]
whence
\beql{1}
J_{s,k}(X)=\int_{(0,1]^k}|f(\bal)|^{2s}d\bal.
\eeq
Our version of the efficient congruencing method will also use the 
exponential sums
\[f_k(\bal;X,\xi,a)=f_a(\bal;\xi)=
\sum_{\substack{x\le X\\ x\equiv\xi\modd{p^a}}}e(\al_1 x+\ldots+\al_k x^k),\]
where $p$ is prime and $a$ is a positive integer exponent.  The prime
$p\ge 5$ will be fixed throughout the argument, so we will not include it
explicitly among the parameters for $f_a(\bal;\xi)$. Taking $s$
and $k$ as fixed we will write
\[I_m(X;\xi,\eta;a,b)=\int_{(0,1]^k}|f_a(\bal;\xi)|^{2m}
|f_b(\bal;\eta)|^{2(s-m)}d\bal,\;\;\;(0\le m\le s-1),\]
which counts solutions of (\ref{eqs}) in which
\[x_i\equiv\xi\modd{p^a}\;\;\; (1\le i\le m\mbox{ and } s+1\le i\le
s+m),\]
and
\[x_i\equiv\eta\modd{p^b}\;\;\; (m+1\le i\le s\mbox{ and } s+m+1
\le i\le 2s).\]
The reader should think of this as a simplified version of Wooley's
$I_{a,b}^{m,r}(X;\xi,\eta)$, given by \cite[(2.11)]{WAp}.  
We observe that when $m=0$ we have
\[I_0(X;\xi,\eta;a,b)=\int_{(0,1]^k}|f_b(\bal;\eta)|^{2s}d\bal,\]
which in independent of $\xi$ and $a$.

We will
also work with $I_m(X;a,b)$ defined by
\[I_0(X;a,b)=\max_{\eta\modd{p}}I_0(X;\xi,\eta;a,b)\]
and
\[I_m(X;a,b)=\max_{\xi\not\equiv\eta\modd{p}}I_m(X;\xi,\eta;a,b)\;\;\;(1\le
m\le s-1). \]
The condition $\xi\not\equiv\eta\modd{p}$ is the last remaining
vestige of Wooley's ``conditioning'' step.  We note for future
reference the trivial symmetry relation
\[I_m(X;a,b)=I_{s-m}(X:b,a).\]

Although many of our results can be proved for general $s$ and $k$ we
shall now specialize to the case $s=6$, $k=3$, and write $J(X)=J_{6,3}(X)$
for brevity.  When $m=0$ we can relate $I_0(X;a,b)$ to $J$ as follows.
\begin{lemma}\label{JI0}
If $p^b\le X$ we have
\[I_0(X;a,b)\le J(2X/p^b).\]
\end{lemma}

We will prove this in the next section along with a number of other
estimates relating different values of $I_1(X;a,b)$ and $I_2(X;a,b)$.
Our next result shows
how to bound $J(X)$ in terms of $I_2(X;1,1)$.
\begin{lemma}\label{JI}
If $p\le X$ we have
\[J(X)\ll pJ(2X/p)+p^{12}I_2(X;1,1).\]
\end{lemma}

One way to compare values of $I_1(X;a,b)$ and $I_2(X;a,b)$
is by applying H\"{o}lder's inequality.  We give two such estimates.
\begin{lemma}\label{hol}
We have
\[I_2(X;a,b)\le I_2(X;b,a)^{1/3}I_1(X;a,b)^{2/3}.\]
\end{lemma}
\begin{lemma}\label{Ik1}
If $p^b\le X$ we have
\[I_1(X;a,b)\le I_2(X;b,a)^{1/4}J(2X/p^b)^{3/4}.\]
\end{lemma}

Next we show how successively larger values of $a$ and $b$ arise.
\begin{lemma}\label{I1k}
We have
\[I_1(X;a,b)\le p^{3b-a}I_1(X;3b,b)\]
if $1\le a\le 3b$.
\end{lemma}

Finally we shall need a result analogous to Lemma \ref{I1k} for $I_2(X;a,b)$.
\begin{lemma}\label{kup}.
If $1\le a\le b$ we have
\[I_2(X;a,b)\le 2bp^{4(b-a)}I_2(X;2b-a,b).\]
\end{lemma}

We are now ready to assemble all these results to prove the following
recursive estimate.
\begin{lemma}\label{rec}
If $1\le a\le b$ and $p^b\le X$ we have
\[I_2(X;a,b)\le 2bp^{-10a/3+14b/3}I_2(X;b,2b-a)^{1/3}
I_2(X;b,3b)^{1/6}J(2X/p^b)^{1/2}.\]
\end{lemma}

For the proof we successively apply Lemmas \ref{kup}, \ref{hol},
\ref{I1k} and \ref{Ik1}, giving
\begin{eqnarray*}
I_2(X;a,b)&\le& 2bp^{4(b-a)}I_2(X;2b-a,b)\\
&\le& 2bp^{4(b-a)}I_2(X;b,2b-a)^{1/3}I_1(X;2b-a,b)^{2/3}\\
&\le& 2bp^{4(b-a)}I_2(X;b,2b-a)^{1/3}\left\{p^{3b-(2b-a)}I_1(X;3b,b)\right\}^{2/3}\\
&\le& 2bp^{4(b-a)+2(a+b)/3}I_2(X;b,2b-a)^{1/3}\\
&& \hspace{1cm}\mbox{}\times
\left\{I_2(X;b,3b)^{1/4}J(2X/p^b)^{3/4}\right\}^{2/3}\\
&=& 2bp^{-10a/3+14b/3}I_2(X;b,2b-a)^{1/3}
I_2(X;b,3b)^{1/6}J(2X/p^b)^{1/2}.
\end{eqnarray*}
Here we should observe that, in applying Lemma \ref{I1k} to
$I_1(X;2b-a,b)$ the necessary condition ``$a\le 3b$'' is satisfied,
since $2b-a\le 3b$.

Everything is now in place to complete the proof of the theorem.  We
note the trivial upper bound $J(X)\ll X^{12}$ and the trivial lower
bound $J(X)\ge [X]^6\gg X^6$ (coming from the obvious diagonal
solutions $x_i=x_{6+i}$ for $i\le 6$).  Thus we may define a real
number $\Delta\in [0,6]$ by setting
\beql{DD}
\Delta=\inf\{\delta\in\R: J(X)\ll X^{6+\delta}\mbox{ for } X\ge
1\}.
\eeq
It follows that we will have $J(X)\ll_{\ep}X^{6+\Delta+\ep}$ for any
$\ep>0$.  Our goal of course is to show that $\Delta=0$.

We observe that
\[I_2(X;a,b)\le J(X)\ll_{\ep} X^{6+\Delta+\ep}\]
for $1\le a\le b$, and hence that
\beql{bc}
I_2(X;a,b)\ll_{\ep} X^{6+\Delta+\ep}p^{-2a-4b}p^{3(3b-a)},
\eeq
since $3(3b-a)\ge 2a+4b$ for $a\le b$.  We now proceed to use Lemma
\ref{rec} to prove, by induction on $n$, that
\beql{ind}
I_2(X;a,b)\ll_{\ep,n,a,b} X^{6+\Delta+\ep}p^{-2a-4b}p^{(3-n\Delta/6)(3b-a)}
\eeq
for any integer $n\ge 0$, provided that
\beql{C1}
1\le a\le b
\eeq
and 
\beql{C2}
p^{3^{n}b}\le X.
\eeq
The base case $n=0$ is exactly the bound (\ref{bc}). 
The reader may be 
puzzled by the choice of the exponent for $p$ in (\ref{ind}).  We shall
discuss this further in the final section.

Given (\ref{ind}) we have
\begin{eqnarray*}
I_2(X;b,2b-a)&\ll_{\ep,n,a,b}&
X^{6+\Delta+\ep}p^{-2b-4(2b-a)}p^{(3-n\Delta/6)(3(2b-a)-b)}\\
&=& X^{6+\Delta+\ep}p^{4a-10b}p^{(3-n\Delta/6)(5b-3a)}.
\end{eqnarray*}
Note that the conditions corresponding to (\ref{C1}) and (\ref{C2}) are
satisfied if 
\[p^{3^{n+1}b}\le X.\]
since we will have $1\le b\le 2b-a$ whenever $1\le
a\le b$, and
\[p^{3^{n}(2b-a)}\le p^{3^{n+1}b}\le X.\]

In a similar way, (\ref{ind}) implies that
\begin{eqnarray*}
I_2(X;b,3b)&\ll_{\ep,n,b}&X^{6+\Delta+\ep}p^{-2b-12b}p^{(3-n\Delta/6)(9b-b)}\\
&=& X^{6+\Delta+\ep}p^{-14b}p^{(3-n\Delta/6)(8b)}
\end{eqnarray*}
the conditions corresponding to (\ref{C1}) and (\ref{C2}) holding
whenever $b\ge 1$.

Finally we have
\[J(2X/p^b)\ll_{\ep} X^{6+\Delta+\ep}p^{-6b-\Delta b}\]
provided that $p^b\le X$. Feeding these estimates into Lemma \ref{rec}
we deduce that
\begin{eqnarray*}
I_2(X;a,b)&\ll_{\ep,n,a,b}& p^{-10a/3+14b/3}
\{X^{6+\Delta+\ep}p^{4a-10b}p^{(3-n\Delta/6)(5b-3a)}\}^{1/3}\\
&& \mbox{}\times 
\{X^{6+\Delta+\ep}p^{-14b}p^{(3-n\Delta/6)(8b)}\}^{1/6}
\{X^{6+\Delta+\ep}p^{-6b-\Delta b}\}^{1/2}\\
&=& X^{6+\Delta+\ep}p^{-2a-4b}p^{(3-n\Delta/6)(3b-a)}p^{-\Delta b/2}\\
&\le &X^{6+\Delta+\ep}p^{-2a-4b}p^{(3-(n+1)\Delta/6)(3b-a)},
\end{eqnarray*}
since $b/2\ge (3b-a)/6$. This provides the required induction step.

Having established (\ref{ind}) we apply it with $a=b=1$, and $p$ chosen 
to lie in the range
\[\tfrac12 X^{1/3^n}\le p\le X^{1/3^n}.\]
There will always be a suitable $p\ge 5$ if 
\[X\ge 10^{3^n}.\]
We then deduce from Lemma \ref{JI}
that
\[J(X)\ll pJ(2X/p)+p^{12}I_2(X;1,1)\ll_{\ep,n} 
p(X/p)^{6+\Delta+\ep}+X^{6+\Delta+\ep}p^{12-n\Delta/3}.\]
If $\Delta$ were strictly positive we could choose $n$ sufficiently large 
that $n\Delta\ge 39$, and would then conclude that
\[J(X)\ll_{\ep,n} X^{6+\Delta+\ep}p^{-1}\ll_{\ep,n}X^{6+\Delta-3^{-n}+\ep},\] 
contradicting the definition (\ref{DD}).  We must therefore have $\Delta=0$,
as required for the theorem.

The reader will probably feel that the final stages oof the argument,
from (\ref{ind}) onward, are lacking in motivation.  The final section
of the paper will offer an explanation for the route chosen.

\section{Proof of the Lemmas}

We begin by examining Lemma \ref{JI0}.  We observe that there is an
$\eta\in (0,p^b]$ such that $I_0(X;a,b)$ counts solutions 
to (\ref{eqs}) in which each
$x_i$ takes the shape $\eta+p^by_i$, with integer variables
$y_i$. We will have $0\le y_i\le X/p^b$.  Thus if we set $z_i=y_i+1$ 
we find that $1\le z_i\le 1+X/p^b\le 2X/p^b$, in view of our condition
$p^b\le X$. Moreover we know that if the $x_i$ satisfy (\ref{eqs}) 
then so to will the $y_i$
and the $z_i$.  It follows that $I_0(X;a,b)\le J_{s,k}(2X/p^b)$ as
claimed.

To prove Lemma \ref{JI} we split solutions of (\ref{eqs}) into
congruence classes for which $x_i\equiv\xi_i\modd{p}$ for $1\le i\le
12$. The number of solutions in which
\[x_1\equiv\ldots\equiv x_{12}\modd{p}\]
is at most 
\[\sum_{\eta\modd{p}}I_0(X;0,\eta;1,1)\le pI_0(X;1,1)\le pJ(2X/p),\]
by Lemma \ref{JI0}. For the remaining solutions to (\ref{eqs})
there is always a pair of variables that are incongruent modulo $p$, and it
follows that there exist $\xi\not\equiv\eta\modd{p}$ such that
\[J(X)\le pJ(2X/p)+
\left(\begin{array}{cc} 12\\ 2 \end{array}\right)p(p-1)
\int_{(0,1]^3}|f_1(\bal;\xi_i)f_1(\bal;\mu)f(\bal)^{10}|d\bal.\]
By H\"{o}lder's inequality we have
\begin{eqnarray*}
\lefteqn{\int_{(0,1]^3}|f_1(\bal;\xi_i)f_1(\bal;\mu)f(\bal)^{10}|d\bal}\\
&\le&
\left\{\int_{(0,1]^3}|f_1(\bal;\xi_i)|^4|f_1(\bal;\mu)|^8|d\bal
\right\}^{1/12}\\
&& \hspace{1cm}\mbox{}\times
\left\{\int_{(0,1]^3}|f_1(\bal;\xi_i)|^8|f_1(\bal;\mu)|^4|d\bal
\right\}^{1/12}\\
&& \hspace{1cm}\mbox{}\times
\left\{\int_{(0,1]^3}|f(\bal)|^{2s}d\bal\right\}^{5/6},
\end{eqnarray*}
whence
\[J(X)\ll pJ(2X/p)+ p^2 I_2(X;1,1)^{1/12}I_2(X;1,1)^{1/12}
J(X)^{5/6}.\]
We deduce that 
\[J(X)\ll pJ(2X/p)+ p^{12}I_2(X;1,1),\]
as required for the lemma.

Lemma \ref{hol} is a trivial application of Holder's inequality.  We
have
\begin{eqnarray*}
I_2(X;\xi,\eta;a,b)&=&\int_{(0,1]^3}|f_a(\bal;\xi)|^4
|f_b(\bal;\eta)|^8d\bal\\
&\le&\left\{\int_{(0,1]^3}|f_a(\bal;\xi)|^8
|f_b(\bal;\eta)|^4d\bal\right\}^{1/3}\\
&& \mbox{}\times
\left\{\int_{(0,1]^3}|f_a(\bal;\xi)|^2
|f_b(\bal;\eta)|^{10} d\bal\right\}^{2/3}\\
&\le& I_2(X;b,a)^{1/3}I_1(X;a,b)^{2/3}\\
&=&I_1(X;a,a),
\end{eqnarray*}
and the lemma follows.

For Lemma \ref{Ik1} we note that
\begin{eqnarray*}
I_1(X;\xi,\eta;a,b)&=&
\int_{(0,1]^3}|f_a(\bal;\xi)|^2|f_b(\bal;\eta)|^{10}d\bal\\
&\le &\left\{
\int_{(0,1]^3}|f_b(\bal;\xi)|^4|f_a(\bal;\eta)|^8d\bal
\right\}^{1/4}\\
&& \hspace{1cm}\mbox{}\times\left\{\int_{(0,1]^3}|f_b(\bal;\eta)|^{12}d\bal
\right\}^{3/4}\\
&\le& I_2(X;b,a)^{1/4}I_0(X;b,b)^{3/4}\\
&\le& I_2(X;b,a)^{1/4}J(2X/p^b)^{3/4},
\end{eqnarray*}
by H\"older's inequality and Lemma \ref{JI0}.

Turning next to Lemma \ref{I1k} we note that $I_1(X;\xi,\eta;a,b)$
counts solutions of (\ref{eqs}) in which $x_i=\xi+p^ay_i$ for $i=1$
and $i=7$, and $x_i=\eta+p^by_i$ for the remaining indices $i$. If
we set $\nu=\xi-\eta$ we deduce that the variables 
\[z_i=\left\{\begin{array}{cc} \nu+p^ay_i, & i=1\mbox{ or }7,\\
p^by_i, & \mbox{otherwise,}\end{array}\right.\]
also satisfy (\ref{eqs}).  In particular, the equation of degree $j=3$
yields
\[(\nu+p^a z_1)^3\equiv(\nu+p^a z_7)^3\modd{p^{3b}}.\]
Now, crucially, we use the fact that $\xi\not\equiv\eta\modd{p}$,
whence $p\nmid\nu$.  It follows that we must have $\nu+p^a z_1
\equiv\nu+p^a z_{s+1}\modd{p^{3b}}$, and hence $z_1\equiv 
z_{s+1}\modd{p^{3b-a}}$.  We therefore have $x_1\equiv x_7\equiv
\xi'\modd{p^{3b}}$ for one of $p^{3b-a}$ possible values of $\xi'$,
so that
\[I_1(X;\xi,\eta;a,b)\le p^{3b-a}I_1(X;3b,b),\]
which suffices for the lemma.

Finally we must handle Lemma \ref{kup}. We note that $I_2(X;\xi,\eta;a,b)$
counts solutions of (\ref{eqs}) in which $x_i=\xi+p^ay_i$ for
$i=1,2,7$ and $8$, and $x_i=\eta+p^by_i$ for the remaining 
indices $i$. As in the proof of Lemma \ref{I1k} we set $\nu=\xi-\eta$
and $z_i=x_1-\eta$, so that the $z_i$ also satisfy (\ref{eqs}).  We
will have $p^b\mid z_i$ for $3\le i\le 6$ and $9\le i\le 12$,
whence
\[(\nu+p^ay_1)^j+(\nu+p^ay_2)^j\equiv(\nu+p^ay_7)^j+(\nu+p^ay_8)^j
\modd{p^{bj}}\;\;\;(1\le j\le 3)\]
with $\nu=\xi-\eta\not\equiv 0\modd{p}$.
We shall use only the congruences for $j=2$ and $3$.  On expanding
these we find that
\beql{c1}
2\nu S_1+p^a S_2\equiv 0\modd{p^{2b-a}}
\eeq
and
\[3\nu^2 S_1+3\nu p^a S_2+p^{2a}S_3\equiv 0\modd{p^{3b-a}},\]
where
\[S_j=y_1^j+y_2^j-y_7^j-y_8^j\;\;\; (j=1,2,3).\]
Eliminating $S_1$ from these yields
\[3\nu p^a S_2 + 2p^{2a}S_3\equiv 0\modd{p^{2b-a}},\]
whence
\[3\nu S_2 + 2p^aS_3\equiv 0\modd{p^{2b-2a}}.\]
Moreover (\ref{c1}) trivially implies that
\[2\nu S_1+p^a S_2\equiv 0\modd{p^{2b-2a}}.\]
It appears that we have wasted some information here, but that turns
out not to be the case.

We now call on the following result, which we shall prove at the end
of this section.
\begin{lemma}\label{count}
With the notations above for $S_j$, let $N(p;a,c)$ denote the number
of solutions $(y_1,y_2,y_7,y_8)$ modulo $p^c$ of the congruences
\[2\nu S_1+p^a S_2\equiv 3\nu S_2 + 2p^aS_3\equiv 0\modd{p^c}.\]
Then if $a\ge 1$ and $c\ge 0$  we will have $N(p;a,c)\le (c+1)p^{2c}$.
\end{lemma}

If $y_i\equiv y_{i0}\modd{p^{2(b-a)}}$ for $i=1,2,7,8$ then $x_i\equiv
  \xi_i\modd{p^{2b-a}}$, with $\xi_i=\xi+p^a y_{i0}$.  The number of
  solutions to (\ref{eqs}) counted by $I_2(X;\xi,\eta;a,b)$ for which
$y_i\equiv y_{i0}\modd{p^{2(b-a)}}$ is then given by
\begin{eqnarray*}
\lefteqn{\int_{(0,1]^3}f_{2b-a}(\bal;\xi_1)f_{2b-a}(\bal;\xi_2)
\overline{f_{2b-a}(\bal;\xi_7)f_{2b-a}(\bal;\xi_8)}
|f_b(\bal;\eta)|^8d\bal}\hspace{1cm}\\
&\le & \int_{(0,1]^3}\left|\prod_{i=1,2,6,7}f_{2b-a}(\bal;\xi_i)\right|
|f_b(\bal;\eta)|^8d\bal\\
&\le &\prod_{i=1,2,6,7}\left\{\int_{(0,1]^3}|f_{2b-a}(\bal;\xi_i)|^4
|f_b(\bal;\eta)|^8d\bal\right\}^{1/4}\\
&\le &\prod_{i=1,2,6,7} I_2(X;\xi_i,\eta;2b-a,a)^{1/4}\\
&\le & I_2(X;2b-a,a),
\end{eqnarray*}
by Holder's inequality.  It then follows from Lemma \ref{count} that
\[I_2(X;a,b)\le N\big(p;a,2(b-a)\big)I_2(X;2b-a,b)\le 
2bp^{4(b-a)}I_2(X;2b-a,a)\]
as required.

It remains to prove Lemma \ref{count}, for which we use induction on $c$.
The base case $c=0$ is trivial. When $c=1$ we have $p\mid S_1$ and
$p\mid S_2$ and the number of solutions is $2p^2-p$, which is also
satisfactory. In general we shall say that a solution
$(y_1,y_2,y_7,y_8)$ is singular if 
\[y_1\equiv y_2\equiv y_7\equiv y_8\modd{p},\]
and nonsingular otherwise. For a nonsingular solution the vectors
\[\nabla(2\nu S_1+p^a S_2)\;\;\;\mbox{and}\;\;\;
\nabla(3\nu S_2+2p^a S_3)\]
are not
proportional modulo $p$, since $a\ge 1$ and $p\nmid 6\nu$.  It follows
that a nonsingular solution $(y_1,y_2,y_7,y_8)$ 
of the congruences modulo $p^c$ will lift to exactly $p^2$ solutions 
modulo $p^{c+1}$.  Thus if we write $N_0(p;a,c)$ for the number of 
nonsingular solutions modulo $p^c$ we will have $N_0(p;a,c)\le
2p^{2c}$, by induction.

For a singular solution we have
\[y_1\equiv y_2\equiv y_7\equiv y_8\equiv\beta\modd{p},\]
say. If we write $y_i=\beta+pu_i$ and
\[S'_j=u_1^j+u_2^j-u_7^j-u_8^j\]
we find that
\[2\nu S_1+p^a S_2=2(\nu+\beta p^a)pS_1'+p^{a+2}S_2'\]
and
\[3\nu S_2+2p^a S_3=
6\beta(\nu+\beta p^a)pS_1'+3(\nu+2\beta p^a)p^2S_2'+2p^{a+3}S_3'.\]
Hence
\[2\nu'pS_1'+p^{a+2}S_2'\equiv 
6\beta\nu'pS_1'+3(\nu'+\beta p^a)p^2S_2'+2p^{a+3}S_3'\equiv 0\modd{p^c}\]
with $\nu'=\nu+\beta p^a\not\equiv 0\modd{p}$. Eliminating $S_1'$ from
the second expression yields
\[3\nu'p^2S_2'+2p^{3+a}S_3'\equiv 0\modd{p^c}\]
and we deduce that
\beql{cc1}
2\nu'S_1'+p^{a+1}S_2'\equiv 0\modd{p^{c-1}}
\eeq
and 
\beql{cc2}
3\nu'S_2'+2p^{a+1}S_3'\equiv 0\modd{p^{c-2}}.
\eeq
Since we are counting values of $y_i$ modulo $p^c$ we have to count
values of $u_i$ modulo $p^{c-1}$.  However any solution of
\[2\nu'S_1'+p^{a+1}S_2'\equiv 3\nu'S_2'+2p^{a+1}S_3'\equiv 0\modd{p^{c-2}}\]
modulo $p^{c-2}$ lifts to exactly $p^3$ solutions of the two
congruences (\ref{cc1}) and (\ref{cc2}) modulo $p^{c-1}$, since
\[\nabla(2\nu'S_1'+p^{a+1}S_2')\equiv 2\nu'(1,1,-1,-1)\not\equiv
0\modd{p}.\]
It follows that (\ref{cc1}) and (\ref{cc2}) have $p^3N(p;a+1,c-2)$
solutions, provided of course that $c\ge 2$ for each of the $p$
possible choices of $\beta$.

We are therefore able to conclude that
\[N(p;a,c)\le N_0(p;a,c)+p^4N(p;a+1,c-2)\le 2p^{2c}+p^4N(p;a+1,c-2)\]
for $c\ge 2$, and the lemma then follows by induction on $c$.

We conclude this section by remarking that in this final inductive
argument, we have estimates of the same order of magnitude for both
the number of singular solutions and the number of nonsingular
solutions. When one tries to generalize the argument to systems of
more congruences the singular solutions can dominate the count in an
unwelcome way.  It is for this reason that Wooley's approach requires a
``conditioning'' step in general, in order to remove singular
solutions at the outset.  Fortunately we just manage to avoid this in
our situation.

\section{Remarks on the Conclusion to the Proof}

This final section is intended to shed some light on the argument that
leads from Lemma \ref{rec} to the theorem.

Suppose one assumes that $J(X)\ll_{\ep} X^{\theta+\ep}$ 
for any $\ep>0$ and that for any
positive integers $a\le b$ one has
\beql{sh}
I_2(X;a,b)\ll_{\ep} X^{\theta+\ep}p^{\alpha a+\beta b}
\eeq
for some constants $\alpha$ and $\beta$, for a suitable range $p\le
X^{\delta(\alpha,\beta)}$, say.

Then Lemma \ref{rec} yields
\[I_2(X;a,b)\ll_b X^{\theta}p^{\alpha' a+\beta' b}\]
for $a\le b$, with new constants 
\[\alpha'=-\tfrac{10}{3}-\tfrac{1}{3}\beta,\;\;\;
\beta'=\tfrac{14}{3}+\tfrac{1}{2}\alpha+\tfrac{7}{6}\beta-\tfrac{1}{2}\theta.\]
We can express this by writing
\[\left(\begin{array}{c}\alpha'\\ \beta'\end{array}\right)=\mathbf{c}+
M\left(\begin{array}{c}\alpha\\ \beta\end{array}\right),\]
with
\[\mathbf{c}=\left(\begin{array}{c}-10/3 \\ 
14/3-\theta/2\end{array}\right)\;\;\;
M=\left(\begin{array}{rr}0 & -1/3\\ 1/2 & 7/6\end{array}\right).\]

Starting with $\alpha=\beta=0$ we can obtain inductively a succession
of bounds of the shape (\ref{sh}), with
\[\left(\begin{array}{c}\alpha\\ \beta\end{array}\right)=
\left(\begin{array}{c}\alpha_n\\ \beta_n \end{array}\right)=\mathbf{c}+
M\mathbf{c}+\ldots +M^n\mathbf{c}.\]
The matrix $M$ has eigenvalues $1$ and $\tfrac{1}{6}$, and can be
diagonalized as $PDP^{-1}$ with
\[P=\left(\begin{array}{rr}-1 & -2\\ 3 & 1 \end{array}\right),\;\;\;
D=\left(\begin{array}{rr}1 & 0\\ 0 & \tfrac{1}{6} \end{array}\right).\]
It then follows that
\[\left(\begin{array}{c}\alpha_n\\ \beta_n \end{array}\right)=
nP\left(\begin{array}{rr}1 & 0\\ 0 &
  0 \end{array}\right)P^{-1}\mathbf{c} +O(1)
=\frac{(6-\theta)n}{5}\left(\begin{array}{c}-1\\ 3\end{array}\right)+O(1)\]
as $n$ tends to infinity. For any starting pair $a,b$ we will have
$3b-a\ge 2b\ge 2$.  Thus if $\theta>6$ we will eventually have
$\alpha_n a+\beta_n b<-1$, say, for suitably large $n$. 

We therefore obtain
\[I_2(X;a,b)\ll_\ep X^{\theta+\ep}p^{-1}\]
for $p\le X^{\delta}$, for some $\delta=\delta_n$ depending on
$\theta$.  This leads to a contradiction, as in \S 2.

We therefore see that the crucial feature of Lemma 7 is that it leads
to a matrix $M$ having its largest eigenvalue equal to 1.  The
corresponding eigenvector is $(\alpha,\beta)=(-1,3)$, and the argument
of \S 2 has therefore been expressed in terms of the linear
combination $3b-a$.

\bigskip
\bigskip

Mathematical Institute,

Radcliffe Observatory Quarter

Woodstock Road

Oxford

OX2 6GG

UK
\bigskip

{\tt rhb@maths.ox.ac.uk}

\end{document}